\documentclass[reqno]{amsart}
\usepackage{amssymb, amsthm, amsmath}

\numberwithin{equation}{section}

\newcommand{\R}{\mathbb{R}} 
\newcommand{\C}{\mathbb{C}}
\newcommand{\Z}{\mathbb{Z}} 
\newcommand{\N}{\mathbb{N}}
 
\newcommand{\Q}{\mathbb{Q}}
\newcommand{\x}{\ma{x}}

\newcommand{\e}{\emph}
\newcommand{\rom}{\mathrm} 
\newcommand{\bfP}{\mathbb{P}}
 
\newcommand{\ma}{\mathbf}
\newcommand{\ben}{\begin{enumerate}}
\newcommand{\een}{\end{enumerate}} 
\newcommand{\eit}{\begin{itemize}}
\newcommand{\beq}{\begin{equation}} 
\newcommand{\eeq}{\end{equation}}
 
\newcommand{\ve}{\varepsilon}
\newcommand{\eps}{\epsilon} 
\newcommand{\mcal}{\mathcal}
\newcommand{\lab}{\label} 
\newcommand{\al}{\alpha}
\newcommand{\D}{\Delta} 
\newcommand{\del}{\delta}

\newcommand{\la}{\lambda}

\newcommand{\h}{\rom{gcd}}
\newtheorem{theorem}{Theorem} 
\newtheorem{lemma}{Lemma}
\newtheorem{pro}{Proposition} 

\newtheorem*{cor*}{Corollary}
 
\renewcommand{\mod}{\hspace{-0.25cm}\pmod}
\newcommand{\tmod}{\hspace{-0.1cm}\pmod}
\newcommand{\colt}[2]{\genfrac{}{}{0pt}{1}{#1}{#2}}
\newcommand{\eqm}[2]{\equiv #1 \mod{#2}}

\newcommand{\tstack}[3]{\colt{#1}{\colt{#2}{#3}}}

\renewcommand{\c}{\ma{c}} 
\newcommand{\y}{\ma{y}}
\newcommand{\z}{\ma{z}} 
\renewcommand{\u}{\ma{u}} 

\renewcommand{\leq}{\leqslant} 
\renewcommand{\geq}{\geqslant}
\renewcommand{\d}{\mathrm{d}}
\newcommand{\Mq}{\|Q\|} 
\newcommand{\mq}{m(Q)}
\newcommand{\CC}{\mathcal{W}} 
\renewcommand{\ss}{\mathfrak{S}(Q)}
\DeclareMathOperator{\supp}{supp}

\theoremstyle{definition}

\begin{document}

\title{Density of integer solutions to diagonal quadratic forms}
\author{T.D. Browning}

\address{School of Mathematics,  
University of Bristol, Bristol BS8 1TW}
\email{t.d.browning@bristol.ac.uk}

\subjclass[2000]{11G35 (11P55, 14G05)}



\begin{abstract}
Let $Q$ be a non-singular diagonal quadratic form in at least four
variables. We provide upper bounds for the number of integer solutions to the equation $Q=0$,
which lie in a box with sides of length $2B$, as $B\rightarrow \infty$. The estimates obtained are completely
uniform in the coefficients of the form, and become sharper as they grow larger in modulus.
\end{abstract}

\maketitle

\section{Introduction}

Let $n \geq 3$ and let $Q \in \Z[x_1,\ldots,x_n]$ be a  non-singular
indefinite quadratic form. Given an arbitrary bounded subset
$\mcal{R}$ of $\R^n$, it is natural to investigate the number of zeros
$\x=(x_1,\ldots,x_n) \in \Z^n$ of the equation $Q(\x)=0$ that are
confined to the region $B \mcal{R}=\{\x\in \R^n: B^{-1}\x\in\mcal{R}\}$, as $B\rightarrow \infty$.  
In the present paper we will focus upon the case of integer solutions to 
diagonal quadratic forms, that lie in the box corresponding to
taking $\mcal{R}$ the unit hypercube in $\R^n$.
Suppose once and for all that
\beq\lab{Q} 
Q(\x)=A_1x_1^2+\cdots+ A_n x_n^2, 
\eeq 
for non-zero integers $A_1,\ldots,A_n$ not all of the same sign, and
write $\D_Q=A_1\cdots A_n$ for the discriminant of $Q$.
Our goal is therefore to  understand the asymptotic behaviour of the
counting function 
$$
N(Q;B)=\#\big\{ \x \in \Z^n: ~Q(\x)=0, ~|\x| \leq B \big\},
$$
where $|\x|=\max_{1\leq i \leq n}|x_i|$ denotes the usual norm on $\R^n$.
Using M\"obius inversion it is then possible
to extract information about the corresponding counting function in
which one is only interested in counting primitive vectors.    This
amounts to counting rational points of bounded height on the quadric
hypersurface $Q=0$ in $\bfP^{n-1}$. 
It will suffice to restrict our attention to primitive quadratic
forms throughout our work, in the sense that $A_1,\ldots,A_n$ 
have greatest common divisor $1$.  

It should come as no surprise that the quantity $N(Q;B)$ has received
substantial attention over the years, to the extent that many authors have
established asymptotic formulae for quantities very similar to
$N(Q;B)$.   Let us define the more general counting function
$$
N_w(Q;B)=\sum w(B^{-1}\x), 
$$
for suitable bounded weight functions $w: \R^n \rightarrow \R_{\geq 0}$
of compact support, where the summation is taken over all $\x \in
\Z^n$ such that $Q(\x)=0$. 
In particular we clearly have $N(Q;B)=N_{w^*}(Q;B)$, where
$w^*=\chi_{[-1,1]^n}$ is the characteristic function of the unit hypercube in
$\R^n$.   Several methods have been developed to study $N_w(Q;B)$ for
appropriate weight functions $w$, and we proceed to discuss what is known.
Under suitable assumptions about $w$, Maly\v sev \cite{malyshev} has
established an asymptotic formula for $N_w(Q;B)$ when $n \geq 5$, and
Siegel \cite{siegel} has done the same when $n=4$ and the discriminant 
$\D_Q$ is a square.  
One of the most impressive results in this direction, however, is due
to Heath-Brown \cite{HB'}, who has established an asymptotic formula for
$N_w(Q;B)$ when $n\geq 3$ and $w$ belongs to a rather general class of
infinitely differentiable weight functions.   Heath-Brown's approach is based upon the
Hardy--Littlewood circle method, and the outcome of his
investigation is the existence of a non-negative constant $c(w;Q)$ such that
\beq\lab{hb}
N_w(Q;B)=c(w;Q)B^{n-2}(\log B)^{b_n}\big(1+o(1)\big),
\eeq
as $B \rightarrow \infty$.  Here $b_n=1$ if $n=3$, or if $n=4$ and $\D_Q$ is a square, and 
$b_n=0$ otherwise.  The constant $c(w;Q)$ may be interpreted as a product of local densities.
All of these estimates for $N_w(Q;B)$ share the common feature that they depend
intimately upon the coefficients of the quadratic form under
consideration.  

The central theme of this paper is the finer question of whether it is
possible to provide estimates for the counting function $N(Q;B)$,
for suitable choices of $Q$, in which the dependence
upon the coefficients of $Q$ is made completely explicit.
Let 
\beq\lab{h} 
\mq=\min_{1\leq i \leq n}|A_{i}|,\quad \Mq=\max_{1\leq i \leq n}|A_{i}|
\eeq 
denote the minimum and height of $Q$, respectively.
For small values of $n$ the geometry of numbers is
particularly effective for this sort of problem.  Thus when $n=3$ it
follows from the author's joint work with Heath-Brown \cite[Corollary
2]{n-2} that 
\beq\lab{one1}
N(Q;B) \ll \Big(\frac{h_Q^{1/2}B }{|\D_Q|^{1/3}} +1\Big) d(|\D_Q|),
\eeq
where $h_Q$ is the greatest common divisor of $A_1A_2,A_1A_3$ and
$A_2A_3$, and $d$ denotes the divisor function.
When $n=4$ work of the author \cite[Theorem~3]{qpub} establishes that
\beq\lab{one2} 
N(Q;B) \ll_\ve \frac{B^{2+\ve}}{\mq^{1/3}|\D_Q|^{1/6}}+ B^{3/2+\ve}, 
\eeq 
for any $\ve >0$, under the assumption that $\D_Q$ is square-free.
With more care, $B^{2+\ve}$ can be
replaced by $B^2 |\D_Q|^{\ve}$ in this estimate.
Both \eqref{one1} and \eqref{one2} have the obvious feature of
becoming sharper as the discriminant of the form grows
larger.  When $n \geq 5$ the best uniform estimate available
is the estimate 
\beq\lab{hb:n-2} 
N(Q;B) \ll_{\ve,n} B^{n-2+\ve}, 
\eeq 
that is due to Heath-Brown \cite[Theorem 2]{annal}.  
Our main result consists of an estimate for $N(Q;B)$ that bridges
\eqref{one2} and \eqref{hb:n-2} for arbitrary $n \geq 4$.

\begin{theorem}\lab{main}
Let $n\geq 4$ and assume that $\D_Q$ is not a square when $n=4$.
Then for any $\ve>0$ and any $B \geq 1$, we have
$$
N(Q;B) \ll_{\ve,n}
\Big(\frac{ B^{n-2}}{\mq^{1/2}\Mq^{1/2}} + 
\frac{\Mq^{2n+3}}{\mq^{3n/4+3}|\D_Q|^{1/2}}
B^{(n-1+\del_n)/2+\ve}\Big)
|\D_Q|^\ve,
$$
where 
\beq\lab{even-odd} 
\del_n=
\left\{\begin{array}{ll} 1, & \mbox{if $n$ is even and $n \geq 5$,}\\ 
0, & \mbox{if $n$ is odd or $n =4$.}
\end{array}
\right.
\eeq 
\end{theorem}

In view of the upper bounds $1\leq \mq$ and $\mq^{-1}|\D_Q|\leq \Mq^{n-1}$, it is clear that we can
always take $\mq^{1/2}\Mq^{1/2} \geq |\D_Q|^{1/(2(n-1))}$ in Theorem \ref{main}.  
For a typical diagonal quadratic form one expects the coefficients to
have equal order of magnitude $|\D_Q|^{1/n}$, so that
there exist positive constants $c_1\geq c_2$, depending
only on $n$, such that
\beq\lab{same-order}
c_1 |\D_Q|^{1/n}\geq  \Mq \geq \mq \geq c_2 |\D_Q|^{1/n}.
\eeq
The following result is a trivial consequence of Theorem
\ref{main}.

\begin{cor*}
Let $n\geq 4$ and assume that $\D_Q$ is not a square when $n=4$.
Suppose that \eqref{same-order} holds for appropriate
constants $c_1,c_2$.  Then for any $\ve>0$ and any $B \geq 1$, we have
$$
N(Q;B) \ll_{\ve,n}
\Big(\frac{ B^{n-2}}{|\D_Q|^{1/n}}
+|\D_Q|^{3/4} B^{(n-1+\del_n)/2+\ve}\Big)
|\D_Q|^\ve,
$$
where $\del_n$ is given by \eqref{even-odd}.
\end{cor*}

A standard probabilistic argument suggests that $N(Q;B)$ should have
order of magnitude $|\D|^{-1/n}B^{n-2}$, at least on average. Our
bounds are clearly consistent with this heuristic. A brief
discussion of certain lower bounds for $N(Q;B)$, in the case $n=4$, can be found 
in the author's  earlier work upon this problem \cite[\S 4]{qpub}.

It is somewhat annoying that the term $\mq$ should appear at all in the
statement of Theorem \ref{main}. We can obtain estimates independent
of $\mq$ by considering an alternative counting function.
Given a parameter $X\geq 1$, let
$$
M(Q;X)=\#\Big\{ \x \in \Z^n: ~Q(\x)=0, ~\max_{1\leq i\leq n}|A_i x_i^2| \leq X \Big\}.
$$
We will deduce the following result rather easily from our proof of
Theorem \ref{main}.

\begin{theorem}\lab{main'}
Let $n\geq 5$.  Then for any $\ve>0$ and any $X \geq 1$, we have
$$
M(Q;X) \ll_{\ve,n}
\Big( \frac{X^{(n-2)/2}}{|\D_Q|^{1/2}} + \Mq^{n/2+\ve}X^{(n-1+\del_n)/4+\ve}\Big) |\D_Q|^\ve,
$$ 
where $\del_n$ is given by \eqref{even-odd}.
\end{theorem}

It would not be hard to extend Theorem \ref{main'} to cover the case
in which $n=4$ and $\D_Q$ is not a square.
Our approach to estimating $N(Q;B)$ and $M(Q;X)$ is based on Heath-Brown's new
version of the Hardy--Littlewood circle method \cite{HB'} that was used to establish
\eqref{hb}.  
Whereas the classical form of the circle method (as
described by Davenport \cite{dav}, for example) is based on the equality 
$$
\int_0^1 e^{2\pi i \al n } \d\al=\left\{
\begin{array}{ll}
1, &\mbox{if $n=0$,}\\
0, &\mbox{if $n\in \Z\setminus\{0\}$,}
\end{array}
\right.
$$
Heath-Brown works with a more sophisticated expression for this
indicator function. The other main difference is the use of Poisson
summation to introduce a family of complete exponential sums, rather
than using the major and minor arc distinction that appears in the classical circle method.

The overall plan will be to establish a version of the asymptotic formula \eqref{hb},
for a suitable weight function $w: \R^n \rightarrow \R_{\geq 0}$, in
which the error terms dependence on $Q$ is made completely explicit.
Once coupled with a uniform upper bound for the constant $c(w;Q)$,
this will suffice for the proof of Theorem \ref{main}.
It is worth highlighting that the classical form of the circle method
could easily be used to establish a result of 
the type in Theorem~\ref{main} when $n\geq 5$. However, a
double Kloosterman refinement is needed to treat the case $n=4$. 
Heath-Brown's approach already incorporates a single
Kloosterman refinement when $n \geq 5$, which in itself yields a sharper
error term. Moreover, the double Kloosterman refinement needed to handle
the case $n=4$ can be carried out with little extra trouble.  
There are a number of extra technical difficulties that need to be dealt with
before Heath-Brown's method can be implemented, however. The most substantial
of these involves pinning down the exact dependence of his estimates for certain
exponential integrals upon the quadratic forms under consideration.

Theorem \ref{main} can be extended in a number of obvious
directions. In addition to  covering the case in which $n=4$ and the
discriminant $\D_Q$ is a square, it is possible to handle non-diagonal indefinite quadratic forms.  
We have decided to pursue neither of these refinements here, however,
choosing instead to focus upon the simplest situation for which we can provide the strongest
results.

We end this section by introducing some of the basic conventions and
notations that we will follow throughout this work.  As is common
practice, we will allow the small positive constant $\ve$ to take
different values at different points of the argument.  We will often
arrive at estimates involving arbitrary parameters $M,N$.  These will
typically be non-negative or positive, but will always take integer values.
Given any vector $\ma{z} \in \R^n$ we write $\int f(\z) \d\z$ for
the $n$-fold repeated integral of $f(\z)$ over $\R^n$.  Given $q \in
\N$, a sum with a condition of the form $\ma{b}\tmod q$ will mean
a sum taken over $\ma{b} \in \Z^n$ such that the components of
$\ma{b}$ run from $0$ to $q-1$. Finally, for any $\al \in \R$ we will
write $e(\al)=e^{2\pi i \al}$ and $e_q(\al)=e^{2\pi i \al/q}$.


\section{Preliminaries}\lab{2}

In this section we bring together the principal ingredients in the
proof of Theorems ~\ref{main} and \ref{main'}. As indicated above, the main idea is to
establish a uniform version of \eqref{hb}, for a suitable weight
function.
Before introducing the weight that we will work
with, we first elaborate upon the nature of
the constant $c(w;Q)$ that appears in Heath-Brown's estimate. As is
well-known to experts, we have $c(w;Q)=\sigma_\infty(w;Q)\ss$, where
$\sigma_\infty(w;Q)$ corresponds to the singular integral, and $\ss$ is the 
singular series. Define the $p$-adic density of solutions to be
\beq\lab{sig-p}
\sigma_p=\lim_{k \rightarrow \infty} p^{-k(n-1)}\#\{\x \mod{p^k}: Q(\x)\eqm{0}{p^k}\},
\eeq
for any prime $p$. When these limits exist, the singular series is given by 
\beq\lab{singseries}
\ss=\prod_{p}\sigma_p.
\eeq
We will see shortly that $\ss$ is convergent for 
the forms considered here.

Consider the function $w_0: \R \rightarrow \R_{\geq 0}$, given by
\beq\lab{w0} 
w_0(x) = \left\{ \begin{array}{ll}  e^{-(1-x^2)^{-1}}, &
|x| < 1,\\  
0, & |x|\geq 1.
\end{array} 
\right.
\eeq 
Then $w_0$ is infinitely differentiable with compact support
$[-1,1]$.  Let
\beq\lab{c0}
c_0=\int_{-\infty}^\infty w_0(x) \d x,
\eeq
and define the function
$$
\omega_\eps(x)=c_0^{-1}\eps^{-1}\int_{-\infty}^{x-\eps}
w_0(\eps^{-1}y)\d y,
$$
for given $\eps >0$.  It is easy to see that $\omega_\eps$ takes values in
$[0,1]$ and is infinitely differentiable, with compact support
$[0,2\eps]$.  In our work we will make use of the non-negative weight function
\beq\lab{dag}
w^{\dag}(\x)=w_0(x_1-2)\prod_{i\neq 1} \omega_{\frac{1}{2}}\Big(1-\frac{x_i}{x_1}\Big),
\eeq
on $\R^n$.  It is clear that $w^{\dag}(\x)$ is zero unless $1\leq x_1 \leq 3$ and
$0 \leq x_i \leq x_1$ for $i \geq 2$.  In particular 
$w^\dag$ is supported in the compact region $[1,3]\times [0,3]^{n-1}$.

We are now ready to record our asymptotic formulae for the weighted counting
function $N_{w^\dag}(Q;B)$.  Let $Q$ be the primitive quadratic form \eqref{Q},
and recall the definitions \eqref{h} of the minimum and height
of $Q$. We may and will assume that the coefficient $A_1$ is positive, throughout 
our work.   The following result will be used to handle the case $n \geq 5$ in Theorem \ref{main}.

\begin{pro}\lab{main-5}
Let $n\geq 5$. Then there exists a non-negative constant
$\sigma_\infty(Q)$ such that 
$$
N_{w^\dag}(Q;B) = \sigma_\infty(Q)\ss B^{n-2}+
O_{\ve,n}\Big( 
\frac{\Mq^{2n+3+\ve}}{A_1^{3n/4+3}|\D_Q|^{1/2}}
B^{(n-1+\del_n)/2+\ve}\Big),
$$
where $\del_n$ is given by \eqref{even-odd}, $\ss$ is given by
\eqref{singseries}, and 
\beq\lab{ox}
\sigma_{\infty}(Q)\ll_n A_1^{-1/2}\Mq^{-1/2}.
\eeq
\end{pro}

Turning to the case $n=4$, for which we must assume that the
discriminant is not a square, we have the the following result.

\begin{pro}\lab{main-4}
Let $n=4$ and assume that $\D_Q$ is not a square.
Then there exists a non-negative constant
$\sigma_\infty(Q)$ such that 
$$
N_{w^\dag}(Q;B) = \sigma_\infty(Q)\ss B^{2}
+O_{\ve}\Big(
\frac{\Mq^{11+\ve}}{A_1^{6}|\D_Q|^{1/2}}
B^{3/2+\ve} \Big),
$$
where $\ss$ is given by
\eqref{singseries}, and $\sigma_{\infty}(Q)$ satisfies \eqref{ox}.
\end{pro}

Our final ingredient in the proof of Theorem \ref{main} is a
uniform upper bound for the singular series $\ss$.  This will show, in
particular, that for the family of quadratic forms \eqref{Q} considered here,
$\ss$ is convergent and actually grows
rather slowly in terms of the coefficients of $Q$. 
The following result will be established in \S \ref{pf:ss-upper}. 

\begin{pro}\lab{ss-upper}
Let $n\geq 4$ and assume that $\D_Q$ is not a square when $n=4$.  Then we have
$
\ss \ll_{\ve,n} |\D_Q|^\ve.
$
\end{pro}

We are now ready to deduce the statement of Theorem \ref{main} from
the statements of Propositions \ref{main-5}--\ref{ss-upper}.
Let $\ve>0$, let $n\geq 4$, and assume that $\D_Q$ is not a
square when $n=4$.  On writing $Q^\sigma$ for the diagonal quadratic form
obtained by permuting the coefficients $A_1,\ldots,A_n$, for each $\sigma
\in S_n$, we deduce that
\begin{align*}
N(Q;B) &\ll_n 1+\sum_{\sigma\in S_n}\sum_{j=0}^\infty N_{w^\dag}(Q^\sigma;B/2^j)\\ 
&\ll_{\ve, n}
\Big(\frac{B^{n-2}}{\mq^{1/2}\Mq^{1/2}} + 
\frac{\Mq^{2n+3+\ve}}{\mq^{3n/4+3}|\D_Q|^{1/2}}
B^{(n-1+\del_n)/2+\ve}\Big)
|\D_Q|^\ve.
\end{align*}
This completes the proof of Theorem \ref{main}.
The proof of Theorem \ref{main'} is handled in exactly the same
way. Instead of using Proposition \ref{main-5}, however, we
employ the main technical result in recent joint work of the author
with Dietmann \cite[Proposition~1]{qleast}. 
Once combined with Proposition \ref{ss-upper}, this latter result
implies that 
$$
N_{w_Q}(Q;B) \ll_{\ve,n} \Big(\frac{B^{n-2}}{|\D_Q|^{1/2}} +
\Mq^{n/2+\ve}B^{(n-1+\del_n)/2+\ve}\Big)|\D_Q|^\ve,
$$
where
$$
w_Q(\x):=w_0(2|A_1|^{1/2}x_1-2)w_0(|A_2|^{1/2}x_2)\cdots w_0(|A_n|^{1/2}x_n),
$$
and $w_0$ is given by \eqref{w0}.  Taking $B=X^{1/2}$, and arguing as
in the deduction of Theorem \ref{main}, we therefore complete the
proof of Theorem \ref{main'}.

It is now time to recall the technical apparatus behind Heath-Brown's
version of the Hardy--Littlewood circle method \cite{HB'}.  Recall the
definitions \eqref{w0} and \eqref{c0} of the
weight function $w_0:\R \rightarrow \R_{\geq 0}$, and the
constant $c_0$.   Let $\omega(x)=4c_0^{-1}w_0(4x-3)$, and define the function $h:
(0,\infty)\times \R \rightarrow \R$ by
$$
h(x,y)=\sum_{j=1}^\infty \frac{1}{xj}\Big(
\omega(xj)-\omega(|y|/xj)\Big).
$$
It is shown in \cite[\S 3]{HB'} that $h(x,y)$ is infinitely
differentiable for $(x,y) \in (0,\infty)\times \R$, and that $h(x,y)$
is non-zero only for $x \leq \max\{1,2|y|\}$.  
Let $Q$ be the quadratic form \eqref{Q}, let $w^\dag$ be given by \eqref{dag}, and let
$X>1$.   The kernel of our work is Heath-Brown's \cite[Theorem
2]{HB'}.   For any $q \in \N$ and any $\c \in \Z^n$, we define the sum
\beq\lab{Sq} S_q(\c) = \sum_{\colt{a=1}{(a,q)=1}}^q
\sum_{\ma{b}\mod{q}}  e_q (aQ(\ma{b})+\ma{b}.\c), 
\eeq 
and the integral 
\beq\lab{Iq} 
I_q(\c)=
\int_{\R^n}w^\dag\Big(\frac{\x}{B}\Big)h\Big(\frac{q}{X},\frac{Q(\x)}{X^2}\Big)e_q(-\c.\x)\d\x.
\eeq 
Then we deduce from the statements of  \cite[Theorems 1 and
2]{HB'} that there exists a positive constant $c_X$, satisfying
$$
c_X=1+O_N(X^{-N}) 
$$
for any integer $N\geq 1$, such that
\beq\lab{asym1} 
N_{w^\dag}(Q;B) = c_X X^{-2}\sum_{\c\in \Z^n} \sum_{q
=1}^\infty q^{-n}S_q(\c)I_q(\c).  
\eeq

In our work we will make the choice 
\beq\lab{X}
X=A_1^{1/2}B, 
\eeq 
where as usual $A_1$ is assumed to be positive.
Things can be made notationally less cumbersome by taking $X=B$ here
instead.  However, this would ultimately lead to a version of
Propositions \ref{main-5} and \ref{main-4} with $A_1$ set to $1$, and
there is no extra technical difficulty in working with \eqref{X}.  In fact
the key property required of $X$ is that we should have
$B^2X^{-2}\partial Q(\x)/\partial x_1\gg_n 1$ on the support of
$w^\dag$.   When $X$ is given by \eqref{X}, we obviously have 
$$
B^2X^{-2}\partial Q(\x)/\partial x_1\geq 2A_1^{-1}A_1 x_1\gg 1 
$$
on $\supp(w^\dag)$,  which is satisfactory. 

Our proof of Propositions \ref{main-5} and \ref{main-4}
now has two major components: the estimation of the exponential sum
\eqref{Sq} and that of the exponential integral \eqref{Iq}.  We will treat the
former in \S \ref{sommes}, while the treatment of the latter is 
rather harder, and will be the focus of \S \ref{integrales}.
We will deduce the statement of Proposition~\ref{main-5} in
\S \ref{eichel}, and that of 
Proposition~\ref{main-4} in \S \ref{eichel'}.
Finally the proof of Proposition \ref{ss-upper} will take place in  \S \ref{pf:ss-upper}.

\section{Estimating $S_q(\c)$}\lab{sommes}

The purpose of this section is to provide good estimates for the exponential sums $S_q(\mathbf{c})$, as given by
\eqref{Sq}, which are uniform in the coefficients of $Q$.
Much of this section follows the general lines of Heath-Brown's investigation \cite[\S\S
9--11]{HB'}.  A number of the results we will need may be 
quoted directly from that work, and we begin by recording the following 
multiplicativity property \cite[Lemma 23]{HB'}.

\begin{lemma}\lab{mult}
If $\h(u,v)=1$ then  
$$ 
S_{uv}(\c)=S_{u}(\c) S_{v}(\c).
$$
\end{lemma}

In fact our work may be further simplified by appealing to the
author's joint work with Dietmann \cite{qleast}, in which uniform estimates for the
average order of $S_q(\mathbf{c})$ are provided for $n \geq 5$.
The outcome of this investigation is the following result
\cite[Lemma 7]{qleast}.

\begin{lemma}\lab{qleast-1}
Let $n \geq 5$ and let $Y\geq 1$.  Then we have
$$
\sum_{q \leq Y}|S_q(\c)| \ll_{\ve,n}
|\D_Q|^{1/2+\ve} Y^{(n+3+\del_n)/2+\ve}, 
$$
where $\del_n$ is given by \eqref{even-odd}.
\end{lemma}

In view of Lemma \ref{mult}, the function $q^{-n}S_q(\ma{0})$ is 
multiplicative.  Moreover, Lemma \ref{qleast-1} implies that the
corresponding infinite sum $\sum_{q=1}^\infty q^{-n}S_q(\ma{0})$ is
absolutely convergent for $n\geq 5$.
Thus the usual analysis of the singular series yields
\beq\lab{yemuna}
\sum_{q=1}^\infty
q^{-n}S_q(\ma{0})=\prod_p \sum_{t=0}^\infty p^{-nt}S_{p^t}(\ma{0})\\
=\prod_p \sigma_p,
\eeq
where $\sigma_p$ is given by \eqref{sig-p}, and we may conclude that 
\beq\lab{singseries1}
\sum_{q \leq Y} q^{-n}S_q(\ma{0}) = \ss +O_{\ve,n}\big(|\D_Q|^{1/2+\ve}Y^{(3+\del_n-n)/2+\ve}\big),
\eeq
for $n\geq 5$. Here, $\ss(Q)$ is given by \eqref{singseries}.

It will suffice to assume that $n=4$ throughout
the remainder of this section.  The following easy upper bound 
for $S_q(\c)$ follows from the proof of \cite[Lemma 25]{HB'}.

\begin{lemma}
\label{expsumm}
We have
$$
S_q(\mathbf{c}) \ll q^{3} \prod_{1 \leq i \leq 4}\h(q,A_i)^{1/2}.
$$
\end{lemma}

\begin{proof}
An application of Cauchy's inequality yields
$$
|S_q(\mathbf{c})|^2 \leq \phi(q)  \sum_{\colt{a=1}{\h(a,q)=1}}^q
\sum_{\ma{d,e}=1}^q e_q \big(
a(Q(\ma{d})-Q(\ma{e}))+\c.(\ma{d}-\ma{e})\big).
$$
On substituting $\ma{d}=\ma{e}+\ma{f}$, we find that
$$
e_q \big( a(Q(\ma{d})-Q(\ma{e}))+\c.(\ma{d}-\ma{e})\big)= e_q (
aQ(\ma{f})+\c.\ma{f})e_q(a\ma{e}.\nabla Q (\ma{f})).
$$
Since $\nabla Q (\ma{f})=2(A_1f_1,\ldots,A_4f_4)$, the summation over
$\ma{e}$ will produce a contribution of zero unless $q \mid A_if_i$
for $1 \leq i \leq 4$. This condition clearly holds for  $\ll
\h(q,A_1)\cdots\h(q,A_4)$ values of $\ma{f} \tmod{q}$,  whence the result.
\end{proof}

We will be able to improve upon Lemma \ref{expsumm} when $q$ is
square-free. The first step is to examine the sum at prime 
values of $q$.  Define the quadratic form 
$$
Q^{-1}(\ma{y}) =A_1^{-1} y_1^2+A_2^{-1} y_2^2+ A_3^{-1} y_3^2+
A_4^{-1} y_4^2,
$$
with coefficients in $\Q$.  When $p$ is a prime such that $p \nmid
2\D_Q$ we may think of $Q^{-1}$ as being defined modulo $p$. 
With this in mind, we have the following result.

\begin{lemma}\lab{s-p}
Let $p$ be an odd prime.  
Then we have
$$
S_p(\c)= \left\{\begin{array}{ll}
-\Big(\frac{\D_Q}{p}\Big)p^{2}, &   \mbox{if $p \nmid Q^{-1}(\c)$},\\
\Big(\frac{\D_Q}{p}\Big)p^{2}(p-1), &   \mbox{if $p \mid Q^{-1}(\c)$},
\end{array}
\right.
$$
if $p\nmid \D_Q$, and 
$$
S_p(\c) \ll p^{5/2}\h(p,\D_Q Q^{-1}(\c))^{1/2}\prod_{i=1}^4\h(p,A_i)^{1/2},
$$
if $p\mid \D_Q$.
\end{lemma}

\begin{proof}
The first part follows on taking $n=4$ in \cite[Lemma 26]{HB'}. The
second part follows by arguing along the
lines of \cite[Lemma 5]{qleast}.
\end{proof}

We may now combine Lemma \ref{mult} and Lemma \ref{s-p} to provide an
estimate for $S_q(\mathbf{c})$ in the case that $q$ is square-free.

\begin{lemma}
\label{expsumm-sfree}
Let $q \in \N$ be square-free. Then we have
$$
S_q(\mathbf{c}) \ll_\ve  q^{5/2+\ve} \h(q,\D_Q Q^{-1}(\c))^{1/2}
\prod_{1 \leq i \leq 4}\h(q,A_i)^{1/2}.
$$
\end{lemma}

\begin{proof}
Since $q$ is square-free we may write $q=2^e\prod_{j=1}^r p_j$, with
$p_1,\ldots,p_r$ distinct odd primes and $e \in \{0,1\}$.  Then it
follows from Lemma \ref{mult}, together with the trivial bound
$|S_2(\c)|\leq 2^4$, that  
$$
|S_q(\c)| \leq 2^{4}
\prod_{j=1}^r |S_{p_j}(\c)|.
$$
Now for each $1 \leq j \leq r$, it follows from Lemma \ref{s-p} that 
$$
S_{p_j}(\c) \ll 
\left\{
\begin{array}{ll}
p_j^{5/2} \h(p_j, \D_Q Q^{-1}(\c))^{1/2}, & \mbox{if $p_j \nmid \D_Q$,}\\
p_j^{5/2} \h(p_j, \D_Q Q^{-1}(\c))^{1/2}\prod_{i=1}^4\h(p_j,A_i)^{1/2}, & \mbox{if $p_j \mid \D_Q$.}
\end{array}
\right.
$$
Putting these estimates together therefore yields the proof of Lemma \ref{expsumm-sfree}.
\end{proof}

We are now ready to discuss the average order of $|S_q(\c)|$, as a
function of $q$. 

\begin{lemma}\lab{exp-av-q}
Let $Y\geq 1$. Then we have
$$
\sum_{q \leq Y}|S_q(\c)| \ll_{\ve}
\left\{
\begin{array}{ll}
|\D_Q|^{1/2+\ve}|\c|^\ve Y^{7/2+\ve}, 
& 
\mbox{if $Q^{-1}(\c)\neq 0$,}\\
|\D_Q|^{1/2} Y^{4}, & 
\mbox{otherwise.}
\end{array}
\right.
$$
\end{lemma}

\begin{proof}
The second bound is an easy consequence of 
Lemma \ref{expsumm}. We therefore proceed under the assumption that
$Q^{-1}(\c)\neq 0$. Write $q=uv$ for
coprime $u$ and $v$, such that $u$ is square-free and $v$ is
square-full.   Then we may combine Lemmas \ref{mult},
\ref{expsumm} and \ref{expsumm-sfree} to deduce that
\begin{align*}
S_q(\c) 
&\ll  |S_{u}(\bar{v}\c)| v^{3} \prod_{1\leq i \leq 4}\h(v,A_i)^{1/2} \\ 
&\ll_\ve u^{5/2+\ve}v^{3} \h(u,\D_Q Q^{-1}(\c))^{1/2} \prod_{1\leq i \leq
4}\h(uv,A_i)^{1/2}\\ 
&\ll_\ve |\D_Q|^{1/2}
q^{5/2+\ve}v^{1/2} \h(u, \D_Q Q^{-1}(\c))^{1/2}.
\end{align*}
It follows that
\begin{align*}
\sum_{q \leq Y}|S_q(\c)| &\ll_\ve  |\D_Q|^{1/2} Y^{5/2+\ve}\sum_{v \leq
Y}  v^{1/2}\sum_{u \leq Y/v} \h(u,\D_Q Q^{-1}(\c))\\
&\ll_\ve |\D_Q|^{1/2+\ve} |\c|^{\ve}
Y^{7/2+\ve}\sum_{v \leq Y}v^{-1/2},
\end{align*}
provided that $Q^{-1}(\c)\neq 0$.
We complete the proof of Lemma \ref{exp-av-q} by noting that there are
$O(V^{1/2})$ square-full values of $v \leq V$. 
\end{proof}

Lemma \ref{exp-av-q} will suffice for our purposes if $Q^{-1}(\c)\neq 0$. 
To handle the case in  which $Q^{-1}(\c)= 0$ we must work somewhat
harder. Consider the Dirichlet series
\beq\lab{goat}
D(s;\c)=\sum_{q=1}^\infty q^{-s}S_q(\c),
\eeq
for $s=\sigma+it \in \C$.
Then it follows from Lemma \ref{exp-av-q} that $D(s;\c)$ is absolutely convergent for $\sigma>4$.  Moreover Lemma 
\ref{mult} yields $D(s;\c)= \prod_p D_p(s;\c)$, where
\beq\lab{msri}
D_p(s;\c)= \sum_{k=0}^\infty p^{-ks}S_{p^k}(\c).
\eeq
We now investigate the factors $D_p(s;\c)$ more carefully, for which
we must distinguish between whether or not $p$ is a divisor
of $2\D_Q$.

Suppose first that $p \mid 2\D_Q$.  Then one easily
deduces from Lemma \ref{expsumm} that 
\begin{align*}
D_p(s;\c)
&\ll \sum_{k=0}^\infty p^{(3-\sigma)k} p^{(\min\{\nu_p(A_1),k\}+\cdots+\min\{\nu_p(A_4),k\})/2}\\
&\ll \sum_{k=0}^\infty p^{(7/2-\sigma)k} p^{\nu_p(\D_Q)/2}
\ll p^{\nu_p(\D_Q)/2},
\end{align*}
if $\sigma>7/2$,  where $\nu_p(n)$ denotes the $p$-adic order of any
non-zero integer $n$. Hence 
\beq\lab{zeta-1}
\prod_{p \mid 2\D_Q} D_p(s;\c) \ll_{\ve} |\D_Q|^{1/2+\ve}.
\eeq
Suppose now that $p \nmid 2\D_Q$.  Then it follows from Lemmas
\ref{expsumm} and \ref{s-p} that
$$
D_p(s;\c)=1+\Big(\frac{\D_Q}{p}\Big)p^{2-s}(p-1) + O_{\delta}(p^{-1-2\del}),
$$
for $\sigma \geq 7/2+\del$, since $Q^{-1}(\c)= 0$.
On writing $\chi_Q(p)=(\frac{\D_Q}{p})$,  we therefore deduce that
\beq\lab{fri=4}
D_p(s;\c)=\big(1-\chi_Q(p)p^{3-s}\big)^{-1} \big(1+ O_{\delta}(p^{-1-\del})\big),
\eeq
for $\sigma \geq 7/2+\del$.  We may combine this with
\eqref{zeta-1} to conclude that
\beq\lab{didas}
D(s;\c)=L(s-3,\chi_Q) E(s;\c),
\eeq
in this region, where $E(s;\c)\ll_{\del,\ve}|\D_Q|^{1/2+\ve}$.
In particular $D(s;\c)$ has an analytic continuation to the half-plane
$\sigma>7/2.$  

Let $Y$ be half an odd integer.  Then it follows from
an application of Perron's formula (see the proof of Titchmarsh
\cite[Lemma 3.12]{titch}, for example), together with the second
estimate in Lemma \ref{exp-av-q}, that 
\beq\lab{zeta-3}
\sum_{q\leq Y} S_q(\c)=\frac{1}{2\pi i} \int_{6-iT}^{6+iT}
\frac{D(s;\c)Y^s}{s}\d s +O\Big(\frac{|\D_Q|^{1/2}Y^{6}}{T}\Big),
\eeq
for any $T \geq 1$.  Let $\al=7/2+\ve$.  Then we proceed to move the line of integration back to
$\sigma=\al$.  Now \eqref{didas} yields
\beq\lab{zeta-2}
D(\sigma+ it;\c)\ll_{\ve} |\D_Q|^{1/2+\ve}|L(\sigma-3 +it,\chi_Q)|,
\eeq
for $\sigma \in [\al,6]$.
We now require the following simple upper bound for the size of the
Dirichlet $L$-function.

\begin{lemma}\lab{crit}
Let $\chi$ be a Dirichlet character modulo $k$. Then we have
$$
L(\sigma+it, \chi) \ll_{\ve, \sigma} 
\left\{
\begin{array}{ll}
k^{(1-\sigma)/2+\ve}\tau^{1-\sigma+\ve}, 
&\mbox{if $\sigma \in (1/2,1]$ and $\chi$ non-principal,}\\
1, 
&\mbox{if $\sigma\in (1,\infty)$,}
\end{array}
\right.
$$
where $\tau=1+|t|$.
\end{lemma}

\begin{proof}
The result is trivial for $\sigma>1$. Assuming that $\sigma\in(1/2,1]$, therefore,
we may combine the P\'olya--Vinogradov inequality with partial
summation, to obtain
\begin{align*}
L(\sigma+it,\chi)
&\ll_\sigma \sum_{n \leq x}\frac{1}{n^\sigma}
+(1+|t|)\int_x^\infty \frac{k^{1/2}\log k}{u^{\sigma+1}}\d u\\
&\ll_\sigma x^{1-\sigma} \log x +\frac{(1+|t|)k^{1/2}\log k}{x^\sigma},
\end{align*}
for any $x\geq 1$.
The proof of the lemma is completed by taking $x=(1+|t|)k^{1/2}$ and
noting that $\log z\ll_\ve z^\ve$ for any $z \geq 1$.
\end{proof}

Sharper versions of Lemma \ref{crit} are available in the literature,
although we will not need anything so deep  here.  For example, 
Heath-Brown \cite{moment} has shown that 
$L(\sigma+it, \chi) \ll_{\ve, \sigma} (k\tau)^{3(1-\sigma)/8+\ve}, $
for $\sigma \in (1/2,1]$, where $\tau=1+|t|$ and $\chi$ is any non-principal character
modulo $k$.

We continue with our analysis of the Dirichlet series $D(s;\ma{c})$.
Now $\chi_Q$ is a non-principal character, since $\D_Q$ is 
not a square when $n=4$.  Hence applying Lemma~\ref{crit} in \eqref{zeta-2} yields
$$
D(\sigma+it;\c)\ll_{\ve}  |\D_Q|^{1/2+\ve}
\left\{
\begin{array}{ll}
|\D_Q|^{(4-\sigma)/2+\ve}(1+|t|)^{4-\sigma+\ve}, & \sigma \in [\al,4],\\
1, & \sigma \in (4,6].
\end{array}
\right.
$$
In particular 
\begin{align*}
\int_{\al \pm iT}^{6 \pm iT}
\frac{D(s;\c)Y^s}{s}\d s 
&\ll_{\ve}  |\D_Q|^{1/2+\ve} \Big( |\D_Q|^{1/4}\frac{Y^{7/2+\ve}
T^\ve}{T^{1/2}} + \frac{Y^6}{T}\Big).
\end{align*}
Turning to the contribution from the vertical lines, we will employ the
mean-value estimate
$$
\int_0^U |L(\sigma+it,\chi)|^2 \d t \ll_\sigma k^{1/2}U,
$$
that is valid for any $\sigma \in (1/2,1)$ and any character modulo
$k$.  But then it follows from this, together with an application of
\eqref{zeta-2} and Cauchy's
inequality,  that
\begin{align*}
\int_{\al - iT}^{\al+ iT}
\frac{D(s;\c)Y^s}{s}\d s 
&\ll_{\ve} |\D_Q|^{1/2+\ve} Y^{7/2+\ve}
\int_{0}^{T} \frac{|L(1/2+\ve+it,\chi_Q)|}{1+t}\d t\\
&\ll_{\ve} |\D_Q|^{3/4+\ve} Y^{7/2+\ve} T^\ve.
\end{align*}
We are now in a position to bring this all together in 
\eqref{zeta-3}.  Thus we conclude the proof of the following result by
taking $T=Y^{5/2}$,  and noting that $D(s;\c)Y^s/s$ is holomorphic in the half-plane $\sigma \geq
\al$.

\begin{lemma}\lab{exp-av-q:4}
Suppose that $\D_Q$ is not a square and $Q^{-1}(\c)= 0$.
Then we have 
$$
\sum_{q \leq Y}S_q(\c) \ll_{\ve} |\D_Q|^{3/4+\ve} Y^{7/2+\ve}.
$$
\end{lemma}

We conclude this section with a few words about the sum $\sum_{q\leq
  Y}q^{-n}S_q(\ma{0})$  in the case $n=4$.  
In the notation of \eqref{goat}, we have $\sum_{q=1}^\infty
q^{-4}S_q(\ma{0})=D(4;\ma{0})$,  and the argument used to prove Lemma
  \ref{exp-av-q:4} ensures that $D(4;\ma{0})$ is convergent. Thus
\eqref{yemuna} continues to hold when $n=4$. Moreover,
we may trace through our application of Perron's formula to conclude that
\beq\lab{sheep}
\sum_{q\leq Y}q^{-4}S_q(\ma{0})=\ss +O_\ve\big(|\D_Q|^{3/4+\ve}
Y^{-1/2+\ve}\big).
\eeq

\section{Estimating $I_q(\c)$}\lab{integrales}

Let $q \in \N$, let $\c \in \Z^n$ and recall the definition \eqref{Q} of the
quadratic form $Q$.  We continue to employ the notation
$\Mq$ for the height of $Q$, as given by \eqref{h}, and the convention
that $A_1>0$ in \eqref{Q}.
The goal of this section is to study
the integral \eqref{Iq}.  In fact it will be convenient to investigate
the behaviour of the integral
\beq\lab{ayl}
I_q(\c;w)=
\int_{\R^n}w\Big(\frac{\x}{B}\Big)h\Big(\frac{q}{X},\frac{Q(\x)}{X^2}\Big)e_q(-\c.\x)\d\x,
\eeq
for a rather general class of weight functions $w:\R^n\rightarrow
\R_{\geq 0}$.  Our first task is therefore to define the class $\CC(S)$ of weight
functions that we will work with. Here, $S$ is an arbitrary set of
parameters that we always assume to contain $n$.

Our presentation will be  much along the lines of \cite[\S\S 2,6]{HB'}.
By a weight function $w$, we will henceforth mean a non-negative function $w: \R^n \rightarrow
\R_{\geq 0}$, which is infinitely differentiable and has compact
support. Given such a function $w$, we let $\dim(w)=n$ denote the
dimension of the domain of $w$ and  
$\rom{Rad}(w)$ be the smallest $R$ such that $w$
is supported in the hypercube $[-R,R]^n$.  Moreover for each integer
$j \geq 0$ we let
$$
\kappa_j(w)=\max\Big\{ \Big|
\frac{\partial^{j_1+\cdots+j_n}w(\x)}{\partial^{j_1}x_1\cdots
\partial^{j_n}x_n}\Big|: ~\x \in \R^n, ~j_1+\cdots+j_n=j\Big\}.
$$
We define $\CC_1(S)$ to be
the set of weight functions $w: \R^n \rightarrow \R_{\geq 0}$, such
that 
$$
\dim(w),\rom{Rad}(w),
\kappa_0(w), \kappa_1(w), \ldots$$
are all bounded by corresponding quantities involving parameters from
the set $S$.  In particular it is clear that $w^\dag \in \CC_1(n)$,
where $w^\dag$ is given by \eqref{dag}.
We now specify the set of functions $\CC(S)\subset \CC_1(S)$. 
Given $w \in \CC_1(S)$, we will say that $w \in \CC(S)$ if 
$x_1 \gg_S 1$ on $\supp(w)$.  In particular it then follows that there is exactly one solution to the 
equation $Q(x_1,\y)=0$, for given $\y=(x_2,\dots,x_n)\in \R^n$ such
that 
$$
\Big(\sqrt{-A_1^{-1}(A_2x_2^2+\cdots+ A_nx_n^2)},\ma{y}\Big) \in \supp(w).
$$  
We conclude our discussion of the class $\CC(S)$ by noting that
$w^\dag\in \CC(n)$.

Returning to the task of estimating $I_q(\c)=I_q(\c;w^\dag)$, it will
clearly suffice to estimate \eqref{ayl} for any choice of weight $w\in
\CC(S)$.  On recalling that $X=A_1^{1/2} B$ in \eqref{X}, a simple change of variables
yields
$$
I_q(\c;w)=
B^n\int_{\R^n}w(\x)h(A_1^{-1/2}B^{-1}q,R(\x))e_q(-B\c.\x)\d\x,
$$
where 
$$
R=A_1^{-1}Q \in \Q[\x].
$$
In particular $R(\x) \ll_S \Mq/A_1$ for any $\x \in \supp(w)$.
It follows from the properties of $h$ discussed in \S \ref{2},
together with the definition of the set $\CC(S)$, 
that $I_q(\c)$ will vanish unless  $q\ll_S B\Mq/A_1^{1/2}$.
Following Heath-Brown we proceed by defining
\beq\lab{I*}
I_r^*(\ma{v};w)= \int_{\R^n}w(\x)h(r,R(\x))e_r(-\ma{v}.\x)\d\x,
\eeq
so that
\beq\lab{I^*}
I_q(\c;w)= B^n I_{r}^*(\ma{v};w),
\eeq
with $r=A_1^{-1/2}B^{-1}q$ and $\ma{v}=A_1^{-1/2}\c$.  We note that
\beq\lab{I^**}
\frac{\partial I_q(\c;w)}{\partial q }= \frac{B^{n-1}}{A_1^{1/2}}
\frac{\partial I_{r}^*(\ma{v};w)}{\partial r}.
\eeq

In our work we will need good upper bounds for 
the integral $I_r^*(\ma{v};w)$ and its first derivative with respect to $r$,
that are uniform in the coefficients of $Q$. 
For this purpose it will clearly suffice to assume that  
$q\ll_S B\Mq/A_1^{1/2}$, or equivalently that $r\ll_S \Mq/A_1$.
Following Heath-Brown's argument  in \cite[\S 7]{HB'},
let $\mcal{H}$ denote the set of infinitely differentiable
functions $f:(0,\infty)\times \R\rightarrow \C$ such that for each $N
\in\N$ there  exist absolute constants $K_{j,N}>0$ for which 
\beq\lab{property}
\Big| \frac{\partial^j f(r,y)}{\partial y^j}\Big| \leq 
\left\{
\begin{array}{ll}
K_{0,N}\big(r^N+\min \big\{1, (r/|y|)^N\big\}\big), & \mbox{if $j=0$,}\\
K_{j,N} r^{-j}\min \big\{1, (r/|y|)^N \big\}, & \mbox{if $j \geq 1$.}
\end{array}
\right.
\eeq
Define the integral
$$
J_r(\u;\omega,f)=\int_{\R^n} \omega(\x)f(r,R(\x))e(-\u.\x)\d\x,
$$
for any $\omega \in \CC(S)$ and any $f \in \mcal{H}$.  Here 
$r$ is restricted to the interval $(0, \infty)$ and $\u$ can be
any vector in $\R^n$.  Let $k \in \{0,1\}$.
Then a straightforward examination of the proof of
\cite[Lemma 14]{HB'} reveals that there exists $\omega^{(k)} \in \CC(S)$ and $f^{(k)} \in
\mcal{H}$ such that $\supp (\omega^{(k)})\subseteq \supp(w)$ and 
\beq\lab{c>0:1:suffice'}
\frac{\partial^k I_r^*(\ma{v};w)}{\partial r^k}
\ll
r^{-1-k}|J_r(r^{-1}\ma{v};\omega^{(k)},f^{(k)})|. 
\eeq
Indeed the only thing to check here is that the statement of \cite[Lemma 14]{HB'}
remains valid when one starts with an arbitrary weight function
belonging to $\CC(S)$, and that it produces auxiliary weight
functions also belonging to $\CC(S)$.

In view of \eqref{c>0:1:suffice'} it will now be enough to estimate
$J_r(\u;\omega,f)$ for given $\omega \in \CC(S)$ and $f \in \mcal{H}$.
We begin by recording a rather trivial upper bound for this integral. The following result is established
much as in \cite[Lemma~15]{HB'}.

\begin{lemma}\lab{J-trivial}
Let $\omega \in \CC(S)$ and let $f \in \mcal{H}$.
Then we have
$
J_r(\u;\omega,f) \ll_S r.
$
\end{lemma}

\begin{proof}
Since $f \in \mcal{H}$, we 
may deduce from \eqref{property} that
$$
J_r(\u;\omega,f) \ll_N \int_{\R^n} \omega(\x)\Big(r^N+ \min\Big\{1, \frac{r^N}{|R(\x)|^N}\Big\}\Big)\d\x,
$$
for any $N\geq 1$.  If $r \geq 1$ then we may take $N=1$ in this
estimate to deduce that $J_r(\u;\omega,f)\ll_S 1+r \ll_S r$, which is
satisfactory for the lemma.  If $r < 1$ then we take $N=2$ to obtain 
$$
J_r(\u;\omega,f) \ll 
\int_{\R^n} \omega(\x)\Big(r+ \min\Big\{1, \frac{r^2}{R(\x)^2}\Big\}\Big)\d\x.
$$
Since $x_1 \gg_S 1$ for any $\x \in \supp(\omega)$ we have 
\beq\lab{lower-partial}
\frac{\partial R}{\partial x_1} \gg_S 1,
\eeq
on $\supp(\omega)$.  On substituting $y=R(\x)$ for $x_1$
we therefore obtain
\begin{align*}
J_r(\u;\omega,f) &\ll_S
r +\int_{|y|\leq r} I(y;\omega)\d y +r^2\int_{|y|> r}
\frac{I(y,\omega)}{y^2}\d y,
\end{align*}
where 
\beq\lab{Iy}
I(y)=I(y;\omega)=\int_{-\infty}^\infty\cdots \int_{-\infty}^\infty
\omega(\x) \frac{\d x_2 \cdots \d x_n}{\partial 
R/\partial x_1},
\eeq
in which $x_1$ is defined by the relation $y=R(\x)$.  
We proceed to show that $I \in \CC_1(S)$, with $\dim(I)=1$.  To see this it suffices to check that
$$
\frac{\partial^j }{\partial y^j}\frac{\omega(x_1(y),x_2,\ldots,x_n)}{\partial
R(x_1(y),x_2,\ldots,x_n)/\partial x_1} \ll_S 1
$$
on $\supp(\omega)$, for any $j \geq 0$.    But this follows from the
lower bound \eqref{lower-partial} and the fact that $\partial
x_1/\partial y= (\partial R/\partial x_1)^{-1}$.  
Having established that $I \in \CC_1(S)$, we obtain
$J_r(\u;\omega,f) \ll_S r$ when $r<1$, as required to complete 
the proof of Lemma \ref{J-trivial}.
\end{proof}

Turning to a more sophisticated treatment of $J_r(\u;\omega,f)$, for
given $\omega \in \CC(S)$ and $f \in \mcal{H}$, we define
\beq\lab{K}
K=A_1^{-1}\Mq, \quad K(\omega)=nK \rom{Rad}(\omega)^2.
\eeq
Then it is clear that $K\ll_S K(\omega) \ll_S K$ for any $\omega \in
\CC(S)$, and 
$
|R(\x)| \leq K(\omega)
$ 
for any $\x \in \supp(\omega)$.  
Recall the definition \eqref{w0} of  $w_0$, and define 
$$
\omega_2(v)=w_0\Big(\frac{v}{2K(\omega)}\Big), \quad \omega_1(\x)=  \frac{\omega(\x)}{\omega_2(R(\x))}.
$$  
Then it is not hard to check that 
$\omega_2$ has compact support $[-2K(\omega),2K(\omega)]$, and that 
$\omega_1 \in \CC(S)$,
with $\supp(\omega_1) \subseteq
\supp(\omega)$.
An examination of the proof of 
\cite[Lemma 17]{HB'} reveals that
\beq\lab{J}
J_r(\u;\omega,f)=\int_{-\infty}^\infty p(t) \int_{\R^n}
\omega_1(\x)e(tR(\x)-\u.\x)\d\x\d t,
\eeq
with 
$$
p(t)=\int_{-\infty}^\infty \omega_2(v)f(r,v)e(-tv)\d v.
$$
We proceed to establish the bound 
\beq\lab{p}
p(t)\ll_{N,S} K r(r|t|)^{-N},
\eeq
for any $N \geq 0$, where $K$ is given by \eqref{K}.
Writing $g(v)=\omega_2(v)f(r,v)$, a repeated 
application of integration by parts reveals that
$$
p(t) \ll_{N,S} |t|^{-N} \int_{-2K(\omega)}^{2K(\omega)} \Big|\frac{d^N
g(v)}{d v^N}  \Big|\d v,
$$
for any $N \geq 0$.  But on employing the inequalities
\eqref{property} satisfied by $f(r,v)$, together with the fact that
$r\ll_S K$, we easily deduce that for any $M\geq 1$ and $N\geq 0$ we have
\begin{align*}
\frac{d^N g(v)}{d v^N} 
&\ll_{M,N,S} \frac{r^M +\min\{1,(r/|v|)^M\} }{K^N} +  
\frac{\min\{1,(r/|v|)^2\}}{r^N}\\
&\ll_{N,S} 
\left\{
\begin{array}{ll}
r^{1-N} +r^{-N}\min\{1,(r/|v|)^2\}, & \mbox{if $r<1$,}\\
r^{1-N}, & \mbox{if $r\geq 1$.}
\end{array}
\right.
\end{align*}
It is now straightforward to deduce the estimate in \eqref{p}.

We are now ready to use the above analysis to deduce a series of
useful basic estimates for the integral $I_q(\c;w)$, for any $w\in \CC(S)$.  Our first result
in this direction will be used to show that large values of $\c$ 
make a negligible contribution in our analysis.

\begin{lemma}\lab{c>0:1}
Let $\c \in \Z^n$ with $\c \neq \ma{0}$.  Then we have
$$
I_q(\c;w)\ll_{N,S}  \frac{B^{n+1}}{q} \frac{\Mq^{N+1}}{A_1^{N/2+1/2}|\c|^{N}},
$$
for any $N\geq 0$.
\end{lemma}

\begin{proof}
The proof of Lemma \ref{c>0:1} closely follows the proof
of \cite[Lemma 19]{HB'}. Recall \eqref{I*}, \eqref{I^*} and the definition \eqref{K} of $K$.  
Then in order to establish Lemma
\ref{c>0:1} it will suffice to show that
\beq\lab{c>0:1:suffice}
I_r^*(\ma{v};w) \ll_{N,S} K^{N+1}r^{-1}|\ma{v}|^{-N},
\eeq
for any $N\geq 0$. To deduce \eqref{c>0:1:suffice}
we employ the identity \eqref{J} and the estimate \eqref{p}
for $p(t)$.
Suppose first that $|\u| \gg_S K|t|$, where $K$ is given by \eqref{K}.  Then we apply \cite[Lemma
10]{HB'} with $f(\x)=tR(\x)-\u.\x$ and $\la=|\u|$. This gives
$$
\int_{\R^n}
\omega_1(\x)e(tR(\x)-\u.\x)\d\x \ll_{M,S} |\u|^{-M},
$$
for any $M\geq 1$. Once inserted into \eqref{J}, and combined with an
application of \eqref{p} with $N=0$, we obtain a contribution of 
$$
\ll_{M,S} K r\int_{|t|\ll_S
K^{-1}|\u|}|\u|^{-M}\d t
\ll_{M,S} r |\u|^{1-M}
$$
to $J_r(\u;\omega,f)$.  When  $|\u| \ll_S K|t|$ we use the trivial bound
$$
\int_{\R^n}
\omega_1(\x)e(tR(\x)-\u.\x)\d\x \ll_{S} 1,
$$
and take $N=M$ in \eqref{p}. This contributes
$$
\ll_{M,S} K r^{1-M}\int_{|t|\gg K^{-1}|\u|}|t|^{-M}\d t
\ll_{M,S} K^M r^{1-M} |\u|^{1-M}
$$
to $J_r(\u;\omega,f)$.
We may now combine these estimates in \eqref{J} to conclude that
$$
J_r(\u;\omega,f) \ll_{M,S}  r|\u|^{1-M} + 
K^{M}r^{1-M}|\u|^{1-M}\ll_{M,S} K^{M}r^{1-M}|\u|^{1-M},
$$
for any $M\geq 1$, since $r \ll_S K$.   Thus it follows that 
$$
J_r(\u;\omega,f) \ll_{N,S} K^{N+1}r^{-N}|\u|^{-N}
$$
for any $N\geq 0$. We may insert this into \eqref{c>0:1:suffice'} with
$k=0$ to deduce that \eqref{c>0:1:suffice} holds for any $N \geq 0$. 
This completes the proof of Lemma~\ref{c>0:1}.
\end{proof}

We will need a finer estimate for $I_q(\c;w)$ when $\c$ has small modulus.  
The following result is established along the lines of \cite[Lemma 22]{HB'}.

\begin{lemma}\lab{c>0:2}
Let $\c \in \Z^n$ with $\c\neq \ma{0}$.
Then we have
$$
q^k \frac{\partial^k I_q(\c;w)}{\partial q^k}\ll_{\ve,S} 
\frac{\Mq^{3n/2+1+\ve}}{A_1^{n/2+1}|\D_Q|}
B^{n/2+1+\ve} \Big(\frac{|\c|}{q}\Big)^{1-n/2+\ve},
$$
for $k \in \{0,1\}$. 
\end{lemma}

\begin{proof}
Our starting point arises from \eqref{I^*}, \eqref{I^**} and \eqref{c>0:1:suffice'}, which render
it sufficient to study the quantity $J_r(\u;\omega,f)$ for $\omega
\in \CC(S)$ and $f\in \mcal{H}$.   We will use the identity \eqref{J}, 
where as stated there $\omega_1 \in \CC(S)$ is such that
$\supp(\omega_1)\subseteq \supp(\omega)$, and $p(t)$ satisfies
\eqref{p} for any $N \geq 0$.

It will be convenient to introduce parameters $\delta>0$ and $T\geq
1$, to be selected in due course.  On recalling the definition \eqref{K} of the
quantity $K$, our immediate goal is to estimate
$J_r(\u;\omega,f)$ under the assumption that
\beq\lab{assume-u}
|\u| \geq KT^2.
\eeq
We proceed by using the subdivision process detailed in \cite[Lemma
2]{HB'} to split up the range for $\x$.  On combining this result with
\eqref{J} it therefore follows that
$$
J_r(\u;\omega,f)=\del^{-n}\int_{-\infty}^\infty p(t) \int_{\R^n}
\int_{\R^n}\omega_\del\Big(\frac{\x-\y}{\del},\y\Big)e(tR(\x)-\u.\x)\d\x\d\y\d t,
$$
where $F(\x)=\omega_\del(\del^{-1}(\x-\y),\y)$ belongs to $\CC(S)$ and
has $\supp(F)\subseteq \supp(\omega_1)\subseteq \supp(\omega)$.
Writing $\x=\y+\del \z$ we deduce that
\beq\lab{J-upper}
|J_r(\u;\omega,f)|\leq \int_{\R^n}\int_{-\infty}^\infty |p(t)| \Big|\int_{\R^n}
\omega_3(\z)e(tR(\x)-\u.\x)\d\z\Big|\d t\d\y,
\eeq
with
$
\omega_3(\z)=\omega_3(\z,\y)=\omega_\del(\z,\y) \in \CC(S)
$
and $\supp(\omega_3) \subseteq \supp (\omega)$.  In order to
effectively estimate the inner integral in \eqref{J-upper}, we must now
differentiate the pairs $(\y,t)$ according to whether or not they
yield a negligible estimate. 
We will say that the pair $(\y,t)$ is `good' if
\beq\lab{good}
\del |t\nabla R(\y)-\u| \geq T \max\{1, K|t|\del^2\}, 
\eeq
and `bad' otherwise.

We begin by treating the case of good pairs $(\y,t)$.  Working under
the assumption \eqref{assume-u} and \eqref{good}, we will estimate the inner integral in \eqref{J-upper}
via an application of \cite[Lemma 10]{HB'}. Let 
$$
f(\z)=tR(\y+\del \z)-\u.\y-\del \u.\z.
$$
Then the partial derivatives of order at least
two are all $O_S(K|t|\del^2)$, and furthermore 
$$
|\nabla f(\z)|=\del |t\nabla R(\y)-\u|
+O_S(K|t|\del^2).
$$
But then $|\nabla f(\z)|\gg_S T \max\{1, K|t|\del^2\}$, and it is easy
to deduce that 
$$
\int_{\R^n}
\omega_3(\z)e(tR(\x)-\u.\x)\d\z \ll_{N,S} T^{-N},
$$
for any $N\geq 1$.  In view of \eqref{p} and the fact
that all of the relevant $\y$ satisfy $|\ma{y}|\ll_S 1$, we therefore
obtain an overall contribution of $O_{N,S} (K T^{-N})$ to
$J_r(\u;\omega,f)$ from the good pairs in \eqref{J-upper}. 

Turning to the contribution from the bad pairs, we henceforth set
\beq\lab{delta}
\del=K^{-1/2}|\u|^{-1/2}.
\eeq
Then it follows from \eqref{good} that 
$$
|t\nabla R(\y)-\u| \leq K^{1/2}T |\u|^{1/2}\max\{1, |t|/|\u|\}, 
$$
if $(\y,t)$ is a bad pair.
We claim that if $(\y,t)$ is a bad pair then 
\beq\lab{order-t}
K^{-1}|\u| \ll_S |t| \ll_S |\u|,
\eeq
and 
\beq\lab{bad}
|t\nabla R(\y)-\u| \ll_S K^{1/2}T |\u|^{1/2}.
\eeq
It will clearly suffice to establish \eqref{order-t}, since
\eqref{bad} is a trivial consequence of this and the previous
inequality. Suppose first that $|t|\leq |\u|$. Then 
$t\nabla R(\y)=\u+O_S(K^{1/2}T |\u|^{1/2})$, and so
$$
K|t| \gg_S |t\nabla R(\y)|\gg_S |\u|,
$$ 
since \eqref{assume-u} implies $|\u|\gg_S K^{1/2}T |\u|^{1/2}$.
This establishes \eqref{order-t} in this case.  Suppose now
that $|t|\geq |\u|$, so that $\u/t=\nabla
R(\y)+O_S(K^{1/2}T|\u|^{-1/2})$. Recall the definition of the weight function 
$$
\omega_3(\z)=\omega_\del(\z,\y)=c_0^{-n} \omega_1(\x)\prod_{i=1}^n w_0\Big(\frac{x_i-y_i}{\del}\Big),
$$ 
as it is constructed in the proof of \cite[Lemma 2]{HB'}.  
In particular it follows that we must have $x_i-\del \leq y_i \leq x_i+\del$ for $1
\leq i \leq n$, if $\omega_3(\z)$ is to be non-zero in \eqref{J-upper}.
Hence $|\nabla R(\y)| \gg_S 1-\del \gg_S 1 \geq K^{1/2}T|\u|^{-1/2},$ by
\eqref{assume-u}.  Thus we must have $|\u|\gg_S |t|$
under the assumption that $T \gg_S 1$, and so \eqref{order-t} holds in
this case also.

Drawing all of this together we deduce from \eqref{J-upper} that 
$$
J_r(\u;\omega,f)\ll_{N,S} K T^{-N}+ \int_{\R^n}\int_{-\infty}^\infty |p(t)| \int_{\R^n}
|\omega_3(\z)| \d\z\d t\d\y,
$$
provided that \eqref{assume-u} holds, 
where $(\y,t)$ runs over values of $\R^{n+1}$ such that
\eqref{order-t} and \eqref{bad} hold, with $|\y|\ll_S 1$.
We now substitute $\x=\y+\del \z$ for $\y$ in this estimate. Now it is
clear from \eqref{delta} that
$$
|t\nabla R(\x) - t\nabla R(\y)| \ll_S K|t|\del \ll_S K^{1/2}|\u|^{1/2}
\leq K^{1/2}T|\u|^{1/2},
$$
since $|\z|\ll_S 1$ and $|t|\ll_S |\u|$.
Thus if $\y$ satisfies \eqref{bad} then so must $\x$.  On employing 
the bound $p(t)\ll_S Kr$ that follows from \eqref{p}, we therefore conclude that
there exists $t\in \R$ in the range \eqref{order-t} such that
$$
J_r(\u;\omega,f)\ll_{N,S} K T^{-N}+  K r|\u|\rom{Vol}(\mcal{S}_t),
$$
for any $N\geq 1$, where
$$
\mcal{S}_t= \{\x \in \supp(\omega): ~ 
|t\nabla R(\x)-\u| \ll_S K^{1/2}T |\u|^{1/2}\}.
$$
An easy calculation reveals that 
each $x_i$ in $\mcal{S}_t$ is restricted to an interval of length
$O_S(|A_i|^{-1}A_1 K^{3/2} T|\u|^{-1/2})$, whence
$$
\rom{Vol}(\mcal{S}_t)\ll_S \frac{K^{3n/2}A_1^{n}}{|\D_Q|} |\u|^{-n/2}T^n.
$$
It therefore follows that
\beq\lab{J-upper'}
J_r(\u;\omega,f)\ll_{N,S} KT^{-N}+
\frac{K^{3n/2+1}A_1^{n}}{|\D_Q|} r|\u|^{1-n/2}T^n,
\eeq
for any $N\geq 1$, under the assumption that \eqref{assume-u} holds.

In order to complete the proof of Lemma \ref{c>0:2} we let $\ve\in
(0,1/2)$ and recall the definition \eqref{K} of $K$.  
We will show that
\beq\lab{J-upper''}
J_r(\u;\omega,f)\ll_{\ve,S} \frac{K^{3n/2+1}A_1^{n}}{|\D_Q|} r|\u|^{1-n/2}(r^{-1}|\u|)^\ve.
\eeq
Suppose first
that $|\u| \leq K^{(n+2\ve)/n}r^{-2\ve/n}$.  Then we have
$$
K^{(n+2\ve)/n}r^{-2\ve/n}\ll_S K^{n/(n-2-2\ve)}r^{-2\ve/(n-2-2\ve)},
$$
since $r \ll_S K$ and $K \geq 1$.
Hence it follows from Lemma \ref{J-trivial} that
$$
J_r(\u;\omega,f)\ll_S r \leq Kr \ll_S K^{n/2+1} r|\u|^{1-n/2}(r^{-1}|\u|)^\ve,
$$
which is satisfactory for \eqref{J-upper''}.  Suppose now that 
$|\u| > K^{(n+2\ve)/n}r^{-2\ve/n}$ and write 
$$
T=c(r^{-1}|\u|)^{\ve/2n},
$$
for a suitable constant $c>0$ depending only on $n$.  Then we
claim that $T\gg_S 1$ and $|\u|\geq KT^2$ if $c$ is chosen to be large enough.
To see the former inequality we note that 
$r^{-1}|\u| > K^{(n+2\ve)/n}r^{-1-2\ve/n}\gg_S 1$, since $r\ll_S
K$. Moreover, the latter inequality holds if and only if 
$|\u| \geq c^{2n/(n-\ve)}K^{n/(n-\ve)} r^{-\ve/(n-\ve)}$.
But this is easily seen to hold when $|\u| >
K^{(n+2\ve)/n}r^{-2\ve/n}$ and $c$ is chosen to be suitably large,
since $r\ll_S K$. Hence we may apply \eqref{J-upper'} to deduce that
$$
J_r(\u;\omega,f)\ll_{N,S} KT^{-N}+
\frac{K^{3n/2+1}A_1^{n}}{|\D_Q|} r|\u|^{1-n/2}(r^{-1}|\u|)^\ve,
$$  
in this case.  On taking $N$ to be sufficiently large in terms of
$\ve$, we therefore complete the proof of \eqref{J-upper''}
for any $\u \in \R^n$.
We now insert this into \eqref{c>0:1:suffice'}, and then into
\eqref{I^*} and \eqref{I^**}, with $r=A_1^{-1/2}B^{-1}q$ and
$\ma{v}=A_1^{-1/2}\c$, in order to conclude the proof of Lemma~\ref{c>0:2}.
\end{proof}

We end this section by considering the integral
$I_q(\ma{0};w)=B^nI_r^*(\ma{0};w)$, for $r=A_1^{-1/2}B^{-1}q$ and any $w \in \CC(S)$.
A trivial application of Lemma \ref{J-trivial}, together
with \eqref{c>0:1:suffice'}, therefore yields
$\partial^k I_r^*(\ma{0};w)/\partial r^k \ll_S
r^{-k}$ for $k \in \{0,1\}$. Hence we may conclude from \eqref{I^*} and \eqref{I^**} that
\beq\lab{I0-triv}
\frac{\partial^k I_q(\ma{0};w)}{\partial q^k} \ll_S q^{-k}B^n,
\eeq 
for $k \in \{0,1\}$.  We can achieve a finer estimate in the case $k=0$.
Since \eqref{lower-partial} holds on $\supp(w)$, we may proceed as in the proof of Lemma
\ref{J-trivial} to conclude that 
$$
I_r^*(\ma{0};w)=\int_{-\infty}^\infty I(y)h(r,y)\d y,
$$
where $I\in \CC_1(S)$ is given by \eqref{Iy}.  
Thus \cite[Lemma~9]{HB'} yields 
$$
I_r^*(\ma{0};w)= I(0) +O_{N,S}(r^N),
$$
for any $N\geq 1$.  In fact \cite[Lemma~9]{HB'} is stated under the
assumption that $r \leq 1$. It is easy to see, however, that the result holds
trivially if  $r\geq 1$, since then $I(0)\ll_S 1 \leq r^N $.  
The integral $I(0)$ is a non-negative constant that is related to the
singular integral, whose value depends 
only upon the weight $w$ and the quadratic form $Q$.  While the
precise value of $I(0)=I(0;w)$ is unimportant for our purposes we will need
the following upper bound.

\begin{lemma}\lab{did}
Let $w \in \CC(S)$. Then we have
$
I(0;w)\ll_S A_1^{1/2}\Mq^{-1/2}.
$
\end{lemma}

\begin{proof}
On relabelling the coefficients
of $Q$, we may assume without loss of generality that 
$$
Q(\x)=A_1x_1^2-\sigma_2A_2x_2^2-\cdots-\sigma_nA_nx_n^2,
$$ 
with $A_1,\ldots,A_n>0$ and $\sigma_2,\ldots,\sigma_n \in \{-1,+1\}$.  It
therefore follows from a simple change of variables that
\begin{align*}
I(0;w)
\leq A_1\int_{-\infty}^\infty\cdots \int_{-\infty}^\infty \frac{
w(A_1^{-1/2}P(\y),A_2^{-1/2}x_2,\ldots,A_n^{-1/2}x_n) \d \y}{|\D_Q|^{1/2}P(\y)},
\end{align*}
where the integral is over all $\y=(x_2,\ldots,x_n)\in \R^{n-1}$,  and 
$$
P(\y)=\sqrt{\sigma_2x_2^2+\cdots+\sigma_nx_n^2}.
$$
But the integrand here vanishes unless $A_1^{1/2}\ll_S P(\ma{y})\ll_S
A_1^{1/2}$.  Hence we have
$$
I(0;w)\ll_S A_1^{1/2}|\D_Q|^{-1/2}\mcal{V}(Q),
$$ 
where 
$\mcal{V}(Q)$ is the volume of
$\y=(x_2,\ldots,x_n)\in \R^{n-1}$ for which $P(\y)\ll_S A_1^{1/2}$ and
$x_i\ll_S A_i^{1/2}$, for $2\leq i \leq n$.
In order to complete the proof of Lemma  \ref{did} it
will therefore
suffice to show that $\mcal{V}(Q)\ll_S |\D_Q|^{1/2}\Mq^{-1/2}$,

This inequality is trivial if $A_1=\Mq$, since 
$\mcal{V}(Q)\ll_S |A_2\cdots A_n|^{1/2}$. Alternatively, on supposing
without loss of generality that $\Mq=A_n$, we fix values of
$x_2,\ldots,x_{n-1}$ and estimate the volume of $x_n \in \R$ such that
$x_n \ll_S A_n^{1/2}$ and $|x_n^2-\del| \ll_S A_1$, where $\del=\del(x_2,\ldots,x_{n-1})$.
The latter inequality implies that this volume is $O_S(A_1^{1/2})$,
which therefore leads to the overall estimate
$$\mcal{V}(Q)\ll_S |A_1\cdots A_{n-1}|^{1/2}=|\D_Q|^{1/2}\Mq^{-1/2},
$$ 
as required.
\end{proof}

Write $I(0;w^\dag)=A_1\sigma_\infty(Q)$, where $w^\dag\in\CC(n)$ is
given by \eqref{dag}. Then on combining Lemma \ref{did} with our arguments above, 
we have therefore established the following result.

\begin{lemma}\lab{c=0}
There exists a non-negative constant $\sigma_\infty(Q)$ such that
$$
I_q(\ma{0})=A_1\sigma_\infty(Q)B^n  +O_{n,N}\big(A_1^{-N/2}q^NB^{n-N}\big),
$$
for any $N\geq 1$, with 
$$
\sigma_\infty(Q)\ll_n A_1^{-1/2}\Mq^{-1/2}.
$$
\end{lemma}

\section{Derivation of Proposition \ref{main-5}} \label{eichel}

In this section we are going to derive Proposition \ref{main-5}. 
Let $n\geq 5$ and recall the choice
$X=A_1^{1/2}B$ that was made in \eqref{X}.  It therefore  follows from
\eqref{asym1} that
\beq\lab{pro1-first}
N_{w^\dag}(Q;B) = \frac{c_B}{A_1 B^2} \sum_{\c\in \Z^n}
\sum_{q =1}^\infty q^{-n}S_q(\c)I_q(\c),
\eeq
where $c_B=1+O_N(B^{-N})$ for any $N\geq 1$, 
$S_q(\c)$ is given by \eqref{Sq}, and $I_q(\c)$ is given by
\eqref{Iq}.  Let $\ve>0$ and $P\geq 1$.
Then it follows from Lemma 
\ref{qleast-1}, together with Lemma \ref{c>0:1} and the fact that
$I_q(\c)=0$ for $q\gg_n B\Mq/A_1^{1/2}$, 
  that the contribution to the right hand
side of \eqref{pro1-first} from  $|\c|> P$ is 
\begin{align*}
&\ll_{n,N} 
\frac{\Mq}{A_1^{3/2} } B^{n-1}
\sum_{q\ll_n B\Mq/A_1^{1/2}} q^{-n-1}|S_q(\c)|  
\sum_{|\c|> P}  \frac{\Mq^{N}}{A_1^{N/2}|\c|^{N}}\\ 
&\ll_{\ve,n,N}
\frac{\Mq^{1+\ve}|\D_Q|^{1/2}}{A_1^{3/2} } 
B^{n-1+\ve}\frac{\Mq^{N}}{A_1^{N/2}P^{N-n}},
\end{align*}
for any $N> n$.  But this is clearly 
\beq\lab{fri=1}
\ll_{\ve,n,M}
\frac{\Mq^{n+1+\ve}|\D_Q|^{1/2}}{A_1^{n/2+3/2} } 
B^{n-1+\ve}\frac{\Mq^{M}}{A_1^{M/2}P^{M}},
\eeq
for any $M\geq 1$.
Turning to the contribution from $1\leq |\c|\leq P$,
we employ Lemma~\ref{c>0:2} to deduce that 
$$
I_q(\c)\ll_{\ve,n}
\frac{\Mq^{3n/2+1+\ve}}{A_1^{n/2+1}|\D_Q|}
B^{n/2+1+\ve} q^{n/2-1}|\c|^{1-n/2+\ve}.
$$
On combining this with Lemma \ref{qleast-1}, we  therefore obtain
\begin{align*}
\sum_{q=1}^\infty q^{-n}S_q(\c)I_q(\c)
&\ll
\max_{Y \ll_n B\Mq/A_1^{1/2}}\sum_{j\leq \log Y}
\sum_{2^{j-1}< q \leq 2^j}q^{-n}|S_q(\c)I_q(\c)| \\
&\ll_{\ve,n}
\frac{\Mq^{3n/2+1+\ve}}{A_1^{n/2+1}|\D_Q|^{1/2}}
B^{n/2+1+\ve} \Big(\frac{B\Mq}{A_1^{1/2}}\Big)^{(1+\del_n)/2}|\c|^{1-n/2+\ve}\\
&\ll_{\ve,n}
\frac{\Mq^{3n/2+2+\ve}}{A_1^{n/2+3/2}|\D_Q|^{1/2}}
B^{(n+3+\del_n)/2+\ve}|\c|^{1-n/2+\ve},
\end{align*}
when $1\leq |\c|\leq P$.
Summing over such values of $\c$ we therefore deduce that
the contribution to the right hand side of \eqref{pro1-first} from
$1\leq |\c|\leq P$ is 
$$
\ll_{\ve,n} \frac{\Mq^{3n/2+2+\ve}}{A_1^{n/2+5/2}|\D_Q|^{1/2}}
B^{(n-1+\del_n)/2+\ve}P^{n/2+1+\ve}.
$$
Once combined with \eqref{fri=1}, we see that
the overall contribution from $\c\neq \ma{0}$ is 
\begin{align*}
\ll_{\ve,n,M}&
\frac{\Mq^{n+1+\ve}|\D_Q|^{1/2} }{A_1^{n/2+3/2}}
B^{(n-1+\del_n)/2+\ve}
\Big(
\frac{\Mq^{M}B^{(n-1)/2}}{A_1^{M/2}P^{M}} + 
\frac{\Mq^{n/2+1}P^{n/2+1+\ve}}{A_1 |\D_Q|}\Big),
\end{align*}
for any $M\geq 1$.  Taking $M=\lceil (n-1)/(2\ve)\rceil$ and 
$
P= |\D_Q|^{1/M}\Mq B^{\ve}/A_1^{1/2},
$ 
we therefore see that there is a contribution of
\beq\lab{water}
\begin{split}
&\ll_{\ve,n}
\frac{\Mq^{n+1+\ve}|\D_Q|^{1/2} }{A_1^{n/2+3/2}}
B^{(n-1+\del_n)/2+\ve}\Big(\frac{1}{|\D_Q|}+
\frac{\Mq^{n+2}}{A_1^{n/4+3/2} |\D_Q|}\Big)\\
&\ll_{\ve,n}
\frac{\Mq^{2n+3+\ve} }{A_1^{3n/4+3}|\D_Q|^{1/2}} B^{(n-1+\del_n)/2+\ve}
\end{split}
\eeq
to the right hand side of \eqref{pro1-first} from
those $\c\neq \ma{0}$.

It remains to handle the contribution from the case $\c=\ma{0}$.
For this it follows from Lemma \ref{qleast-1} and \eqref{I0-triv} that
for any $Y \geq 1$, we have
\begin{align*}
\sum_{Y/2< q \leq Y}q^{-n}S_q(\ma{0})I_q(\ma{0}) 
&\ll_{\ve,n}
|\D_Q|^{1/2+\ve} B^n Y^{(3+\del_n-n)/2+\ve}.
\end{align*}
On summing over dyadic intervals for $Y$ such that
 $B \leq Y \ll_n B \Mq/A_1^{1/2}$, we deduce that the
overall contribution to the right hand side of \eqref{pro1-first} from
$\c= \ma{0}$  and $q\geq B$, is 
\beq\lab{london}
\ll_{\ve,n} A_1^{-1}|\D_Q|^{1/2+\ve} B^{(n-1+\del_n)/2+\ve}.
\eeq
Finally we note that an application of 
\eqref{singseries1}, together with  Lemma \ref{c=0}, reveals that the contribution 
from $\c= \ma{0}$  and $q\leq B$ is 
\begin{align*}
\frac{c_B}{A_1 B^2}\sum_{q \leq B} q^{-n}S_q(\ma{0})I_q(\ma{0})
=&
\sigma_\infty(Q)\ss B^{n-2}\\
&\quad +O_{\ve,n}(A_1^{-1}|\D_Q|^{1/2+\ve}B^{(n-1+\del_n)/2+\ve})\\
&\quad + O_{n,N}\Big(\frac{B^{n-2-N}}{A_1^{N/2+1}} 
\sum_{q \leq B} q^{N-n}|S_q(\ma{0})|\Big),
\end{align*}
for any $N\geq 1$.    
On selecting $N=(n-3-\del_n)/2$, and
applying Lemma \ref{qleast-1}, we deduce that the error terms in this estimate are
also bounded by \eqref{london}.  
Observe that
$$
\frac{|\D_Q|^{1/2}}{A_1}=\frac{|\D_Q|}{A_1 |\D_Q|^{1/2}}\leq 
\frac{\Mq^n}{A_1 |\D_Q|^{1/2}}\leq 
\frac{\Mq^{2n+3}}{A_1^{3n/4+3}|\D_Q|^{1/2}}.
$$
We may now combine these inequalities with \eqref{water} and \eqref{london} in
\eqref{pro1-first}, in order to complete the proof of Proposition \ref{main-5}.

\section{Derivation of Proposition \ref{main-4}} \label{eichel'}

We proceed as in the previous section, much of which carries
over to this setting. Let  $n=4$ and assume that
$\D_Q$ is not a square.  Then the above argument suffices to handle the terms with 
$|\c| \geq P$, or with $1\leq |\c| \leq P$ and $Q^{-1}(\c)\neq 0$, where
$P=|\D_Q|^{\ve}\Mq B^\ve/A_1^{1/2}$.  Hence it follows that
$$
N_{w^\dag}(Q;B)=\frac{c_B}{A_1 B^2} \sum_{\tstack{\c\in \Z^4
}{Q^{-1}(\c)=0}{|\ma{c}|\leq P}}
\sum_{q \ll B\Mq/A_1^{1/2}} 
\hspace{-0.2cm}
q^{-4}S_q(\c)I_q(\c)
+O_{\ve}\Big(\frac{\Mq^{11+\ve}}{A_1^{6}|\D_Q|^{1/2}}  B^{3/2+\ve}\Big),
$$
where $c_B=1+O_N(B^{-N})$.  
When $1\leq |\c|\leq P$ and $Q^{-1}(\c)=0$ we may use partial
summation, based on Lemmas \ref{exp-av-q:4} and \ref{c>0:2}.  Now the
latter result implies that 
$$
\frac{\partial^k I_q(\c)}{\partial q^k}\ll_{\ve}
\frac{\Mq^{7+\ve}}{A_1^{3}|\D_Q|} B^{3+\ve}q^{1-k} |\c|^{-1+\ve},
$$
for $k \in \{0,1\}$. Once combined with the former, we see that
\begin{align*}
\sum_{Y/2<q \leq Y} q^{-4}S_q(\c)I_q(\c)
&\ll_{\ve}
\frac{\Mq^{7+\ve}}{A_1^{3}|\D_Q|^{1/4}} \frac{B^{3+\ve}Y^{1/2+\ve}}{
|\c|^{1-\ve}}
\ll_{\ve} \frac{\Mq^{15/2+\ve}}{A_1^{13/4}|\D_Q|^{1/4}} \frac{B^{7/2+\ve}}{|\c|^{1-\ve}},
\end{align*}
for any $Y \ll_{\ve} B\Mq/A_1^{1/2}$. Now it follows from \eqref{hb:n-2} that there
are at most $O_{\ve,n}(P^{2+\ve})$ vectors $\c \in \Z^4$ for which $|\c|\leq
P$ and $Q^{-1}(\c)=0$.  Hence we conclude that
\begin{align*}
N_{w^\dag}(Q;B) =& \frac{1}{A_1 B^2} 
\sum_{q \ll B\Mq/A_1^{1/2}} 
q^{-4}S_q(\ma{0})I_q(\ma{0})
+O_{\ve}\Big(\frac{\Mq^{11+\ve}}{A_1^{6}|\D_Q|^{1/2}} B^{3/2+\ve} \Big).
\end{align*}

To handle the contribution from large $q$ we employ partial summation
again, this time based on Lemma \ref{exp-av-q:4} and \eqref{I0-triv}.
Thus we obtain
$$
\sum_{Y/2< q \leq Y}q^{-4}S_q(\ma{0})I_q(\ma{0}) 
\ll_{\ve}  |\D_Q|^{3/4+\ve} B^4 Y^{-1/2+\ve}.
$$
for any $Y \geq 1$.  On summing over dyadic intervals for $Y$ in the
interval $B^{1-\ve} \leq Y \ll B \Mq/A_1^{1/2}$, we easily deduce that terms
with $\c=\ma{0}$ and $q\geq B^{1-\ve}$ contribute $O_{\ve}(A_1^{-1}|\D_Q|^{3/4+\ve}
B^{3/2+\ve})$ to  $N_{w^\dag}(Q;B)$, which is satisfactory.
Finally, for $q \leq B^{1-\ve}$ we may apply 
Lemma \ref{c=0} and \eqref{sheep}, together with the second part of
Lemma \ref{exp-av-q},  in order to conclude that
$$
N_{w^\dag}(Q;B) = \sigma_\infty(Q)\ss B^{2}
+O_{\ve}\Big(\frac{\Mq^{11+\ve}}{A_1^{6}|\D_Q|^{1/2}}  B^{3/2+\ve}\Big).
$$
This completes the proof of Proposition \ref{main-4}.

\section{The singular series} \label{pf:ss-upper}

In this section we establish Proposition \ref{ss-upper}.   Let $n \geq
4$ and assume that $\D_Q$ is not a square when $n=4$.  Then \eqref{yemuna}
holds, and we have
$$
\ss=\prod_p D_p(n;\ma{0}),
$$
in the notation of \eqref{msri}.  
We begin by handling the factors $D_p(n;\ma{0})$, for which $p
\nmid \D_Q$.  Suppose first that $p=2$.  Then an application of
\cite[Lemma 4]{qleast} reveals that  
$$
D_2(n;\ma{0})= 1+O_{\ve,n}\Big(\sum_{k\geq 2}2^{k(1-n/2+\ve)}
\Big)\ll_n 1,
$$
since $n \geq 4$.  Suppose now that $p>2$ and 
write $\chi_Q(p)=(\frac{\D_Q}{p})$, as usual. We may therefore combine \eqref{fri=4}
with the proof of \cite[Eqn. (5.2)]{qleast}, in order to deduce that
$$
D_p(n;\ma{0})=
\left\{
\begin{array}{ll}
(1-\chi_Q(p)p^{-1})^{-1} (1+O(p^{-3/2})), & \mbox{if $n=4$},\\
1+ O_{\ve,n}(p^{-3/2+\ve}), &\mbox{if $n\geq 5$}.
\end{array}
\right.
$$
Now Lemma \ref{crit} implies that $L(1,\chi_Q)\ll_\ve |\D_Q|^\ve$.  Hence
\beq\lab{sheep'}
D_2(n;\ma{0})\prod_{p\nmid 2\D_Q} D_p(n;\ma{0})  \ll_{\ve,n} 
\left\{
\begin{array}{ll}
|\D_Q|^\ve, & \mbox{if $n= 4$ and $\D_Q$ not square},\\
1, & \mbox{if $n\geq 5.$}
\end{array}
\right.
\eeq

We now turn to an upper bound for the factors $\sigma_p=D_p(n;\ma{0})$, for odd $p \mid \D_Q$.
Recall that
$$
\sigma_p=\lim_{k \rightarrow \infty} p^{-k(n-1)}N_k(p), \quad N_k(p)=\#\{\x \mod{p^k}: Q(\x)\eqm{0}{p^k}\}.
$$
After relabelling the indices we may assume that there exists $d\in
\N$ such that
$$
(A_1,\ldots,A_n)=(a_1,\ldots,a_r,p^d b_1,\ldots,p^d b_s),
$$
with $p \nmid a_1\cdots a_r$. Moreover, we may suppose that $r+s=n$
and there is an index $1\leq i \leq s$ such that $p \nmid b_i$.
We have $r,s \geq 1$, since $p \mid \D_Q$ and the
highest common factor of $A_1,\ldots,A_n$ is assumed to be $1$.
Observe that for any fixed integers $a,b$ such that $p \nmid a$, and any $k \in
\N$, the number of positive integers $n\leq p^k$
such that $an^2\equiv b \tmod{p^k}$ is at most $2$.
We proceed to show that
\beq\lab{stag}
N_k(p)\leq 4p^{k(n-1)},
\eeq
for any $k \geq 1$.  If $k \leq d$ then it easily follows that
\begin{align*}
N_k(p) &= p^{ks}\#\{x_1,\ldots,x_r \mod{p^k}: \mbox{$\sum_{i=1}^r
  a_ix_i^2\eqm{0}{p^k}$}\}
\leq 2p^{k(n-1)},
\end{align*}
which is satisfactory for \eqref{stag}.  Assume now that $k >d$.  
On writing $\x=(\y,\z)$ for
$\y=(y_1,\ldots,y_r)$ and $\z=(z_1,\ldots,z_s)$ modulo $p^{k}$, we see that
\begin{align*}
N_k(p) &= \#\{\y,\z \mod{p^k}: \mbox{$\sum_{i=1}^r a_iy_i^2+p^d\sum_{i=1}^s b_iz_i^2 \eqm{0}{p^{k}}$}\}\\
&= \sum_{\colt{\y \mod{p^{k}}}{\exists w \in \Z: ~\sum_{i=1}^r
    a_iy_i^2=p^d w}}
\#\{\z \mod{p^{k}}: \mbox{$\sum_{i=1}^s b_iz_i^2 \eqm{w}{p^{k-d}}$}\}.
\end{align*}
Now the summand here is plainly equal to 
$$
p^{ds}\#\{\z \mod{p^{k-d}}: \mbox{$\sum_{i=1}^s b_iz_i^2
  \eqm{w}{p^{k-d}}$}\}
\leq 2p^{ds+(k-d)(s-1)},
$$
since there exists at least one value of $b_1,\ldots,b_s$ that is not
divisible by $p$.  But then it follows that for $k>d$ we have
\begin{align*}
N_k(p) &\leq 2p^{k(s-1)+d}
\#\{\y \mod{p^{k}}: \mbox{$\sum_{i=1}^r  a_iy_i^2 \eqm{0}{p^d}$}\}\\
&\leq 4 p^{k(s-1)+d} p^{(k-d)r+d(r-1)}=
4 p^{k(n-1)}.
\end{align*}
This too is satisfactory for \eqref{stag}.  We have
therefore shown that $\sigma_p\leq 4$ when $n \geq 4$ and $p\mid
\D_Q$ is an odd prime. Once combined with \eqref{sheep'}, this shows that
$$
\ss\ll_{\ve,n} |\D_Q|^\ve \prod_{p \mid \D_Q} 4  \ll_{\ve,n} |\D_Q|^\ve,
$$
and so completes the proof of Proposition \ref{ss-upper}.

\end{document}